\documentclass[reqno, 11pt]{amsart}
\usepackage{amssymb,amsmath,amsfonts}
\usepackage{amssymb}
\usepackage{amsmath,amssymb}
\usepackage[latin1]{inputenc}
\newtheorem{thm}{Theorem}[section]
\newtheorem{lem}[thm]{Lemma}
\newtheorem{cor}[thm]{Corollary}
\newtheorem{prop}[thm]{Proposition}

\newtheorem{rmk}[thm]{Remark}
\newtheorem{obs}[thm]{Remark}
\numberwithin{equation}{section}
\usepackage{enumerate}
\usepackage{graphicx,color}

\everymath{\displaystyle}
\newcommand{\R}{\mathbb{R}}
\newcommand{\Rn}{\mathbb{R}^{n}}
\newcommand{\N}{\mathbb{N}}

\newcommand{\til}{\widetilde}

\begin{document}

\title[Global existence for a nonlinear system with fractional Laplacian ]{Global existence for a nonlinear system with fractional Laplacian in Banach spaces}

\author{ Miguel loayza and Paulo F. S. Silva}

\address{  Departamento de Matem\'atica, Universidade Federal de Pernambuco - UFPE, 50740-540, Recife, PE,                                                                                Brazil}
\email{miguel@dmat.ufpe.br}

\address{ Departamento de Matem\'atica, Universidade Federal do Rio Grande do Norte - UFRN, 50078-970, Natal, PE,                                                                                Brazil}
\email{paulorfss@ccet.ufrn.br}



\date{\today}

\begin{abstract}
We consider the cauchy problem for the fractional power dissipative equation $u_t+(-\Delta )^{\beta/2} u=F(u)$, where $\beta>0$ and $F(u)=B(u, ...,u)$ and $B$ is a multilinear form on a Banach space $E$. We show a global existence result assuming some properties of scaling degree of the multilinear form and the norm of the space $E$. We extend the ideas used for the treating of the equation to determine the global existence for the system $u_t+(-\Delta)^{\beta/2}= F(v)$, $v_t+(-\Delta )^{\beta/2}v=G(u)$ where $F(u)=B_1(u,...,u), G(v)=B_2(v,...,v)$

\end{abstract}

\maketitle

\section{Introduction}
Let $n \in \N$, $\beta>0$. We consider the Cauchy problem for the semilinear fractional power equation
\begin{equation}\label{In.uno}
    \left\{
    \begin{array}{l}
    \partial_tu+(-\Delta)^{\beta/2}u = F(u), \mbox{ in } \Rn \times (0,\infty)\\
    u(0)=u_0 \mbox{ in } \Rn,
    \end{array}
    \right.
\end{equation}
where $F(u)= B(u,...,u)$ and $B: E^p \to \R$ is a $p-$linear form, where $E$ is a Banach space, $p>1$ and $E^p= E\times ...\times E( p \mbox{ times})$. We also assume that the $p-$linear form $B$ acts on $u$ only with respect to the spacial variable. 


The problem (\ref{In.uno}) models several classical problems, for example 

\noindent
(1) The semilinear fractional power dissipative equation  
$$u_t+(-\Delta)^{\beta/2} u=  \nu u^p.$$

\noindent
(2) The  generalized Hamilton-Jacobi equation
$$u_t+(-\Delta)^{\beta/2}=\nabla u \cdot \nabla u.$$
When $\beta=2$, we have the Hamilton-Jacobi equation.

\noindent
(3) The generalized Navier-Stokes equation
$$u_t+(-\Delta u)^{\beta/2} -(u \cdot \nabla) u+ \nabla P=0, \ \nabla \cdot u=0.$$

\noindent
(4) The generalized convection-diffusion equation
$$u_t+(-\Delta u)^{\beta/2}=a \cdot \nabla (u^p), \ a \in \R^n, a \neq 0$$ 

The case $\beta=2$ for  the semilinear dissipative equation correspond to the semilinear heat equation and has been studied extensively, see for instance \cite{fujita}, \cite{BC},\cite{cazenaveweissler}, \cite{snoussi4}, \cite{wei1}, \cite{wei2}. For $\beta \neq 2$ see for example \cite{miao}, \cite{wuyuan}, \cite{wu}. For the nonlinear Hamilton-Jacobi equation the local well posedness in Lebesgue spaces has been discussed in \cite{Ben}. Concerning the generalized Navier-Stokes equation see \cite{Can}. Global well posedness including self-similar solutions and the large time behavior have been proved for convection-diffusion in \cite{AEZ} and \cite{EZ}.

Local and global existence and large time behavior for solutions of problem (\ref{In.uno}), with $\beta=2$ and $B$ a bilinear form on $E\times E$, were studied in a general context in \cite{karch}. Specifically,  it is assumed that 
\begin{enumerate}[(i)]
\item The norm  $\Vert \cdot \Vert_{E}$ has scaling degree equal to $a$, that is 
\begin{equation}\label{In.esc}
\Vert u_\lambda \Vert_{E}=\lambda^{a}\Vert u \Vert_{E}
\end{equation}
for each $u \in E, \lambda>0$ such that $u_\lambda \in E$, where $u_\lambda(x)=u(\lambda x)$ for $x \in \R^n.$ 

\item The bilinear form $B$ has the following scaling property
$$B((u_1)_\lambda,(u_2)_\lambda)=\lambda^b[B(u_1, u_2)]_\lambda$$
for some $b \in \R, \lambda>0$, and for every $u_i \in E$ so that $(u_i)_\lambda \in E~(i=1,2)$.

\item The Banach spaces $E$ is adequate to problem (\ref{In.uno}), that is, 
\begin{enumerate}[(a)]
\item $\mathcal{S} \subset E \subset \mathcal{S}'$ with continuous injections. 

\item The norm is translations invariant, that is, $\Vert u(\cdot + x)\Vert_E=\Vert u \Vert_E$ for all $u \in E$ and $x \in \R^n$.

\item For all $u, v \in E$, $B(u,v) \in \mathcal{S}'$ and  there exist $T_0>0$ and a function $w \in L^1(0,T_0)$, $w>0$ such that $$\Vert S(t) B(u,v)\Vert \leq w(t)\Vert u\Vert \Vert v \Vert,$$ 
for any $u,v \in E$. Here, $(S(t))_{t\geq 0}$ is the heat semigroup.
\end{enumerate} 
\end{enumerate}

Henceforth, $\mathcal{S}$ denotes the space of Schwartz rapidly decreasing functions and $\mathcal{S}'$ its dual space, that is, the space of tempered distributions. Some examples of adequate spaces are the Lebesgue space $L^p(\R^n)$, the Marcinkiewicz weak $L^p(\R^n)$  space, the Lorentz space $L^{p,q}(\R^n)$ and  the Morrey space $\mathcal{M}^p(\R^n)$.

With these concepts it was established the following local existence result.
\begin{thm}[\cite{karch}] Let $\beta=2$, $E$ a Banach space and let $B$ be a bilinear form on $E \times E$ with scaling degree $\sigma <2$. Let $r>n(2-\sigma)$ and $0\leq \alpha<\min\{1, 2-\sigma-n/r, \sigma+n/p\}$. Suppose that $E$ has the following properties:
\begin{enumerate}[(i)]
\item $E$ is adequate to problem (\ref{In.uno});

\item The norm $\Vert \cdot \Vert_{E}$ has a scaling degree $-\frac{n}{r}$;

\item  $S(t): E \to L^q(\R^n)$ is a bounded operator for every $t>0$ and some $q \in [1, \infty]$

\end{enumerate} 

Let $u_0 \in BE^\alpha$. There exists $T=T(u_0)$ and a unique local in time solution of problem (\ref{In.uno}) on $[0,T)$, which is unique in the space
\begin{equation}\label{In.cla}
\mathcal{F}([0,T), BE^\alpha) \cap \{u:(0, T) \to E; \sup_{0<t<T}t^{\alpha/2}\Vert u(t)\Vert_{E}<\infty \}.
\end{equation}
\end{thm}

The space $\mathcal{F}([0, T], BE^\alpha)=\{u \in L^\infty((0,T), BE^\alpha); u(t)\to u(0) \mbox{ as } t \to 0 \mbox{ in }\mathcal{S}'\}.$ The spaces $BE^\alpha$ is given by 
\begin{equation}\label{In.tre}
BE^\alpha=\{f\in \mathcal{S}'; \Vert f \Vert_{BE^\alpha}= \sup_{t>0} t^{\alpha/2}\Vert S(t)f \Vert_{E}<\infty\}.
\end{equation}

For the global existence we have.
 
\begin{thm}[\cite{karch}] Let $\beta=2$ and $B$ be a bilinear form with scaling degree $\sigma <2$ on $E\times E$. Let $r>0$ satisfying $1-\frac{n}{r}<\sigma<2-\frac{n}{r}$. Denote $\alpha=2-\sigma-\frac{n}{r}$. 
Then, there exists $\epsilon>0$ such that for each $u_0 \in BE^\alpha$ satisfying $\Vert u_0 \Vert _{BE^\alpha}<\epsilon$, there exists a global solution for problem  (\ref{In.uno}) in the space (\ref{In.cla})($T=\infty$).
Moreover, this is the unique solution satisfying $\sup_{t>0}t^{\alpha/2} \Vert u(t)\Vert_{E}< 2 \epsilon.$
\end{thm}

The main objective of this work is extend Karch'result for problem (\ref{In.uno}) with $\beta>0$ and $B$ a $p-$linear form on $E^p$. We recall that the fractional Laplacian $(-\Delta)^{\beta/2}$ is defined by
    \begin{equation}\label{deflf}
    (-\Delta)^{\beta/2}u = \mathcal{F}^{-1}\left(|\xi|^{\beta}\mathcal{F}u \right).
    \end{equation}
where $\mathcal{F}$ and $\mathcal{F}^{-1}$ denote the Fourier transform and its inverse, and they are given by $$\mathcal{F}u(\xi) = \int_{\Rn} e^{-2\pi i \langle x,\xi\rangle}u(x)dx,$$ $$\mathcal{F}^{-1}u(x) = \int_{\Rn} e^{2\pi i \langle x,\xi\rangle}u(\xi)d\xi. $$
    
The semigroup generated by the operator  $(-\Delta)^{\beta/2}$ is defined by
    \begin{equation}\label{dsg0}
    S_{\beta}(t)u=K_{\beta}(t)*u,
    \end{equation}
where $u \in \mathcal{S}'$ and $K_{\beta}(t,\cdot)=\mathcal{F}^{-1}(e^{-t|\xi|^{\beta}}) \in \mathcal{S},~t>0$. Note that for $\beta=2$ we have the well known heat kernel, $K(t,x)=(4\pi t)^{-\frac{n}{2}}e^{-\frac{|x|^2}{4t}}$. For this reason,  $(-\Delta)^{\beta/2}$ can be consider as the generalization of the Laplacian operator $-\Delta$. It is easy see that  $S_{\beta}(t)$ linear and $S_{\beta}(t+s) = S_{\beta}(t)S_{\beta}(s).$

If $BE^\alpha$ is a Banach space and $u_0\in BE^\alpha$, we say that $u \in L^\infty((0,T), BE^\alpha)$ is a solution for problem (\ref{In.uno}) if $u$ verifies, in some sense, the equation
\begin{equation}\label{In.cua}
u(t)=S_\beta(t) u_0 + \int_0^t S_\beta(t-\sigma) B(u(\sigma), ..., u(\sigma))d\sigma
\end{equation}
for every $t \in (0,T)$.  

We begin with the local existence for problem (\ref{In.uno}).
\begin{thm}\label{teo.ex} Let $B$ a $p-$linear($p>1$) form  with scaling degree $\sigma$ and $E \in \mathcal{X}$ an adequate Banach space  to problem  (\ref{In.uno}) with  scaling degree $a$. Assume that $\beta+(p-1)a>\sigma$, and
$$0<\alpha<\min\left \{\beta/p, (\beta-\sigma)/(p-1)+a, -(p-1)a+\sigma\right \}.$$
Let $u_0 \in BE^{\alpha}$. There exit $T>0$ and a unique function $u$ in the space
$$\mathcal{F}([0,T], BE^\alpha)\cap \{u:(0, T)\to E; \sup_{0<t<T} t^{\alpha/\beta}\Vert u(t)\Vert_E<\infty\}.$$
\end{thm}

Our result about global existence is de following.
\begin{thm}\label{teorema1} Let $B$ a $p-$linear($p>1$) form  with scaling degree$\sigma$ and $E \in \mathcal{X}$ an adequate Banach space  to problem  (\ref{In.uno}) with  scaling degree $a$, where
    \begin{equation}\label{alpha}
    \frac{\beta-\sigma}{p-1} - \frac{\beta}{p}<-a<\frac{\beta-\sigma}{p-1}\mbox{ and }\alpha=\frac{\beta-\sigma}{p-1}+a.
    \end{equation}
Let $M,R>0$ such that $R+pKM^{p-1}<M$ with $K>0$ given explicitly by  (\ref{Val.k}). Then for every $u_0 \in BE^{\alpha}$ with $\|u_0\|_{BE^{\alpha}} \leq R$, there exists a unique global solution $u$ of problem (\ref{In.uno}) satisfying
      \begin{equation}\label{BE}
      \Vert u\Vert:=\sup_{t>0}\|u(t)\|_{BE^{\alpha}} + \sup_{t>0}t^{\frac{\alpha}{\beta}}\|u(t)\|_{E}\leq M.
      \end{equation}
Moreover, If  $v$ is other global solution for problem (\ref{In.cua}) satisfying (\ref{BE}) and with initial data $v_0$ with $\|v_0\|_{BE^{\alpha}}\leq R$, then
$$\Vert u-v\Vert \leq [1-p(2M)^{p-1}K]^{-1}\|u_0-v_0\|_{BE^{\alpha}}.$$
In addition, if $pK_1M^{p-1}<1$, where $K_1$ is given by (\ref{Val.k1}), we have $$\lim_{t \infty}t^{\alpha/2}\Vert u(t)-v(t)\Vert_E=0$$if and only if 
 $\lim_{t \to \infty}t^{\alpha/2}\Vert S_\beta(u_0-v_0)\Vert_E=0$
  
\end{thm}

\begin{rmk} Here are some comment concerning Theorem \ref{teorema1}
\begin{enumerate}[(i)]
\item  It is clear that if $\beta=p=2$, Theorem \ref{teorema1} reduces to Karch's result. 

\end{enumerate}
\end{rmk}

We use our arguments to analyze  the semilinear fractional power system 
\begin{equation}\label{In.sis}
    \left\{
    \begin{array}{l}
    \partial_tu+(-\Delta)^{\beta/2}u = B_1(v,..., v), \mbox{ in } \Rn \times (0,\infty)\\
    \partial_tv+(-\Delta)^{\beta/2}v = B_2(u,...,u),\mbox{ in } \Rn \times (0,\infty)\\
    u(0)=u_0, v(0)=v_0 \mbox{ in } \Rn,
    \end{array}
    \right.
    \end{equation}
where the $B_1, B_2$ are multi-linear forms defined on Banach spaces. Problem (\ref{In.sis}) for $\beta=2$ has been considered by various authors, see for example, \cite{escobedoherrero1}, \cite{snoussi2}.

As in the problem (\ref{In.uno}), if $u_0 \in BE^{\alpha_1}, v_0 \in BE^{\alpha_2}$ and $BE^{\alpha_1}, BF^{\alpha_2}$ are Banach spaces, we say that $(u,v) \in L^{\infty}((0,T); BE^{\alpha_1})\times L^{\infty}((0,T); BF^{\alpha_2})$ is a solution of the problem (\ref{In.sis}) if verifies, in some  sense, the following system
    \begin{equation}\label{sistint1}
    \left\{
    \begin{array}{ll}
    u(t) =& \displaystyle S_{\beta}(t)u_0 + \int_0^t S_{\beta}(t-\tau)B_1(v,...,v)d\tau, \\
    v(t) =& \displaystyle S_{\beta}(t)v_0 + \int_0^t S_{\beta}(t-\tau)B_2(u, ..., u)d\tau,
    \end{array}
    \right.
    \end{equation}
for every $t \in (0,T)$.

On the global existence for problem (\ref{In.sis}) we have the following result.
\begin{thm}\label{teo2} Let $E, F$ be Banach spaces , $B_1: E^q \to \R$ and $B_2: F^p \to \R$ two forms with scaling degree $\sigma_1$ and $\sigma_2$ respectively. Assume that $E, F\in \mathcal{X}$ have scaling degree $a$ and $b$ respectively and that $E\times F$ is adequate to system (\ref{In.sis}). Let $pq>1$,
\begin{equation}\label{defalpha12}
    \alpha_1=\frac{\beta(1+q)}{pq-1}+a -\frac{\sigma_1+q\sigma_2}{pq-1}, \  \alpha_2=\frac{\beta(1+p)}{pq-1}+b-\frac{\sigma_2+p\sigma_1}{pq-1}.
    \end{equation}
Suppose that
\begin{enumerate}[(i)]
\item $\alpha_1, \alpha_2>0$.    

\item $\alpha_1+q b<a+\sigma_1$ and $\alpha_2+pa<b+\sigma_2$.

\item $\alpha_1<q\alpha_2$, $\alpha_2<p\alpha_1.$

\end{enumerate}
Let $M,R>0$ so that $R + qM^{q}K_1 + pM^{p}K_2 < M.$ Then, for any $~\Phi=(u_0,v_0) \in BE^{\alpha_1}\times BF^{\alpha_2}$ verifying
$$
\mathcal{N}(\Phi):=\|u_0\|_{BE^{\alpha_1}} + \|v_0\|_{BF^{\alpha_2}} \leq R,
$$
there exists an unique solution $U=(u,v)$ for system (\ref{sistint1}) such that
$$
\Vert U\Vert:=\sup_{t>0}t^{\frac{\alpha_1}{\beta}}\|u(t)\|_{E}+\sup_{t>0}\|u(t)\|_{BE^{\alpha_1}} + \sup_{t>0}t^{\frac{\alpha_2}{\beta}}\|v(t)\|_{F} + \sup_{t>0}\|v(t)\|_{BF^{\alpha_2}}\leq M.
$$
Moreover, if $\overline{U}=(\overline{u},\overline{v})$  is a solution for problem  (\ref{sistint1})  verifying $\Vert \overline U\Vert \leq M$ with initial data $\overline{\Phi}=(\overline{u}_0,\overline{v}_0)$ which verify $\mathcal{N}(\overline \Phi)\leq R$, then
    \begin{equation}\nonumber
    \Vert U-\overline{U}|\Vert \leq [1-\left(qM^{q-1}K_1 + pM^{p-1}K_2\right)^{-1}\mathcal{N}(\Phi-\overline{\Phi}).
    \end{equation}
\end{thm}

\begin{rmk} Here are some comment on Theorem \ref{teo2}.
\begin{enumerate}[(i)]
\item If $\sigma_1=\sigma_2=\sigma$ and $b=[(p+1)/(q+1)]a $, then $\alpha_1=[(\beta-\sigma)(q+1)]/(pq-1)$, $\alpha_2=[(p+1)/(q+1)]\alpha_1$ and conditions (i)-(iv) are reduced to 
$$(\beta-\sigma)\frac{q+1}{pq-1}-\frac{\beta}{p}\gamma< -a<(\beta-\sigma)\frac{q+1}{pq-1},$$
where $\gamma=\min\{1, [p(q+1)]/[q(p+1)]\}$. In particular, for $p=q$ we have condition (\ref{alpha}).

\item Local existence for problem (\ref{In.sis}) can be obtained  modifying slightly the proof of Theorem \ref{teo2}. To do this, we assume that 
\begin{enumerate}
\item $\alpha_1, \alpha_2>0$, $\alpha_1-q\alpha_2+ \beta+qb>a+\sigma_1$, $-p\alpha_1+\alpha_2+\beta+pa>b+\sigma_2.$

\item $\beta+qb>a+\sigma_1$, $\beta+pa>b+\sigma_2$. 

\item $\beta>q\alpha_2$, $\beta>p\alpha_1$.

\item $a+\sigma_1>\alpha_1+qb$ and $b+ \sigma_2>\alpha_2+pa.$
\end{enumerate}
\end{enumerate}
\end{rmk}

\section{Preliminary results }

In this section we extend for  $\beta\neq 2$  the definitions considered by Karch in \cite{karch,karch2} for $\beta=2$. The arguments used in the proof of Propositions \ref{ers}, \ref{BEacompleto} and \ref{esb1} are similar to Karch's arguments, but since we are considering situations in that $\beta$ can be different to two  we give them for completeness. 

\subsection{ Scaling properties} Let  $(E, \Vert \cdot \Vert_E)$ be a Banach space which can be imbedded  continuously in $\mathcal{S}'$. We say that the norm $\Vert \cdot \Vert_E$ has scaling degree equal to $a$ if equality (\ref{In.esc}) holds.  What follows are some examples of Banach spaces with its respective scaling degrees:
\begin{enumerate}[(i)]
     \item Lebesgue spaces, $L^p(\Rn)$: $~\|u_{\lambda}\|_{L^p} = \lambda^{-\frac{n}{p}}\|u\|_{L^p},~1\leq p \leq +\infty$,
     \item Lorentz spaces, $L^{(p,q)}(\Rn)$: $~\|u_{\lambda}\|_{L^{(p,q)}} = \lambda^{-\frac{n}{p}}\|u\|_{L^{(p,q)}},~1\leq p,q \leq +\infty$,
     \item Morrey Homogeneous spaces, $\dot{\mathcal{L}}_q^p(\Rn)$: $~\|u_{\lambda}\|_{\dot{\mathcal{L}}_q^p} = \lambda^{-\frac{n}{p}}\|u\|_{\dot{\mathcal{L}}_q^p},~1\leq q \leq p \leq +\infty$,

\item Besov Homogeneous spaces, $B_{p,q}^s(\Rn)$: $~\|u_{\lambda}\|_{\dot{B}^s_{(p,q)}} = \lambda^{s-\frac{n}{p}}\|u\|_{\dot{B}^s_{(p,q)}},~1\leq p,q\leq +\infty~$ and  $~s<0$
    \end{enumerate}
For a study of these spaces see \cite{bergh} and \cite{grafakos}.

We say that a $p-$linear form $B$ defined on $E^p=E \times E \times ...\times E(p \mbox{ times})$ has scaling degree equal to
$\sigma$ if
\begin{equation}\label{Pre.uno}
B((u_1)_{\lambda},...,(u_p)_{\lambda}) = \lambda^{\sigma}\left[ B(u_1,..,u_p) \right]_{\lambda},
\end{equation}
for every $u_1,u_2,...,u_p \in E$, $~\lambda>0~$ and $~(u_i)_{\lambda}=u_i(\lambda \cdot),~i=1,2,...,p$.

Some examples of $p-linear$ forms, with its respective scaling degree $\sigma$,  are given below:
\begin{enumerate}[(i)]
\item $B(u_1,...,u_p)=u_1\cdot ...\cdot u_p,$ has scaling degree $\sigma=0$. 

\item $B(u_1, ..., u_p)=( u_1 \cdot ...\cdot u_{p-1})  \ {\bf a}\cdot \nabla u_p$ with ${\bf a}\in \R^n$ has scaling degree $\sigma=1$.

\item $B(u_1,u_2)=\nabla u_1 \cdot \nabla u_2$ has scaling degree $\sigma=2.$

\end{enumerate}
In the following result we establish some scaling relations. 
\begin{prop}\label{sso} Let $\lambda>0$ and  $u \in \mathcal{S}'$.  Then, $(-\Delta)^{\beta/2}u_{\lambda}=\lambda^{\beta}[(-\Delta)^{\beta/2}u]_{\lambda}$, $S_{\beta}(t)u_{\lambda}=[S_{\beta}(\lambda^{\beta}t)u]_{\lambda}$ and $t^{\frac{n}{\beta}}K_{\beta}(t,t^{\frac{1}{\beta}}x) = K_{\beta}(1,x)$ for $t>0$, $x \in \R^n$.
\end{prop}

\noindent {\bf Proof. } By the change of variable $\xi'=t^{1/\beta}\xi$ we have
$$
\begin{array}{ll}
    K_{\beta}(t,t^{\frac{1}{\beta}}x) &= \int_{\R^n}e^{2\pi i \langle t^{\frac{1}{\beta}}x,\xi\rangle}e^{-t|\xi|^{\beta}}d\xi \\
    &= t^{-\frac{n}{\beta}}\int_{\R^n}e^{2\pi i \langle x,\xi'\rangle}e^{-|\xi'|^{\beta}}d\xi'.
    \end{array}
$$
By (\ref{deflf}) and the change of variables $y'=\lambda y$ and $\xi'=\lambda^{-1}\xi$ we conclude 
$$ 
\begin{array}{lll}
    [(-\Delta)^{\beta/2}u_{\lambda}](x) &=& \displaystyle \int_{\R^n}e^{2\pi i\langle x,\xi\rangle}|\xi|^{\beta}\mathcal{F}u_{\lambda}(\xi)d\xi \nonumber\\
    &=& \displaystyle \int_{\R^n}e^{2\pi i\langle x,\xi\rangle}|\xi|^{\beta}\left[\int_{\R^n}e^{-2\pi i\langle \xi,y\rangle}u(\lambda y)dy\right]d\xi \nonumber\\
    &=& \displaystyle \lambda^{\beta}\int_{\R^n}e^{2\pi i\langle\lambda x,\xi'\rangle}|\xi'|^{\beta}\left[\int_{\R^n}e^{-2\pi i\langle \xi', y'\rangle}u(y')\lambda^{-n}dy'\right]\lambda^{n}d\xi' \nonumber\\
    &=& \lambda^{\beta}[(-\Delta)^{\beta/2}u_{\lambda}](\lambda x) \nonumber\\
    &=& \lambda^{\beta}[(-\Delta)^{\beta/2}u]_{\lambda}(x) \nonumber\\
    & & \nonumber\\
\end{array}
$$
From (\ref{dsg0})
$$
\begin{array}{lll}
    \left[S_{\beta}(t)u_{\lambda}\right](x) &=& \int_{\R^n}\mathcal{F}^{-1}(e^{-t|\xi|^{\beta}})(x-y)u(\lambda y)dy  \nonumber\\
    &=& \int_{\R^n}\mathcal{F}^{-1}(e^{-t|\xi|^{\beta}})(x-\lambda^{-1}y')u(y')\lambda^{-n}dy'  \nonumber\\
    &=& \int_{\R^n}\left[\int_{\R^n}e^{2\pi i\langle x-\lambda^{-1}y', \xi \rangle}e^{-t|\xi|^{\beta}}d\xi\right]u(y')\lambda^{-n}dy' \nonumber\\
    &=& \int_{\R^n}\left[\int_{\R^n}e^{2\pi i\langle\lambda x-y', \xi'\rangle}e^{-t\lambda^{\beta}|\xi'|^{\beta}}\lambda^nd\xi'\right]\varphi(y')\lambda^{-n}dy' \nonumber\\
    &=& \int_{\R^n}\left[\int_{\R^n}e^{2\pi i\langle \lambda x-y', \xi' \rangle}e^{-t\lambda^{\beta}|\xi'|^{\beta}}d\xi'\right]u(y')dy'\nonumber\\
    &=& [\mathcal{F}^{-1}(e^{-\lambda^{\beta}t|\xi|^{\beta}})*u](\lambda x) \nonumber\\
    &=& [S_{\beta}(\lambda^{\beta}t)u]_{\lambda}(x) \nonumber\\
    & & \nonumber \\
\end{array}
$$

When $u$ is an homogeneous function we have the following result.

\begin{cor} If $E$ is a Banach space with scaling degree equal to $a$ and $u$ is an homogeneous  function with degree equal to $\theta$,  then $S_{\beta}(t)u=\lambda^{-\theta}[S_{\beta}(\lambda^{\beta}t)u]_{\lambda}$ and $\|S_{\beta}(1)u\|_{E} = t^{(a-\theta)/\beta}\|S_{\beta}(t)u\|_{E}$.
\end{cor}

\begin{prop}\label{ers} Let $E_a$ and $E_b$ be Banach spaces with scaling degree $a$ and $b$ respectively. Assume that there exists $t_0>0$ such that $S_{\beta}(t_0): E_a\to E_b$ is bounded, where $S_{\beta}(t)$ is defined by (\ref{dsg0}). Then there exists  $C=C(a,b,\beta, t_0)>0$ so that
    \begin{equation}\nonumber
    \|S_{\beta}(t)  u\|_{E_b}\leq Ct^{\frac{a-b}{\beta}}\|u\|_{E_a}.
    \end{equation}
for all $t>0$. Moreover, if $E_b\subset \mathcal{S}'$ is imbedded continuously and $b<a$, then $E_a=\{0\}$.
\end{prop}

\noindent
{ \bf Proof.} Since $S_\beta(t_0)$ is bounded, there exists  $C'>0$ such that $\|S_{\beta}(t_0)u\|_{E_b}\leq C'\|u\|_{E_a},$
for every $u \in E_a$. By Proposition \ref{sso}, we have
$$
    \lambda^{b}\|S_{\beta}(\lambda t_0)u\|_{E_b}=\|[S_{\beta}(\lambda^{\beta}t_0)u]_{\lambda}\|_{E_b}=\|S_{\beta}(t_0)u_{\lambda}\|_{E_b}\leq C'\|u_{\lambda}\|_{E_a}=C'\lambda^{a}\|u\|_{E_a}.
$$
Thus, $\|S_{\beta}(\lambda^{\beta}t_0)u\|_{E_b}\leq C'\lambda^{a-b}\|u\|_{E_a}$ for every $\lambda>0$.
Setting $\lambda=(t/t_0)^{\frac{1}{\beta}}$ we conclude that $\|S_{\beta}(t)  u\|_{E_b}\leq C t^{\frac{a-b}{\beta}}\|u\|_{E_a}$ for every $t>0$ with $C=C't_0^{(b-a)/\beta}.$

Note that if $b<a$ and $u\in E_a$, then $S_{\beta}(t)  u\to 0$ as $t\to 0$ in $E_b$. Since $E_b$ is imbedded continuously in $\mathcal{S}'$ follows that $S_{\beta}(t)  u\to 0$ as $t\to 0$ in $\mathcal{S}'$. Since $S_{\beta}(t)  u\to u$ in $\mathcal{S}'$ as $t \to 0$, by uniqueness, we conclude that $u=0$.  

\begin{prop}\label{lemaqmq0}
Let $E$ be a Banach space with scaling degree $a$. Suppose that for some $t_0>0$ and $r_0\in [1,+\infty]$ the operator $S_{\beta}(t_0):E \to L^{r_0}(\R^n)$ is bounded. Then  $S_{\beta}(t):E \to L^r(\R^n)$  for all $t>0$ and $r \geq r_0$. Moreover, there exists a constant $C>0$ which does not depend of $t>0$ and
    \begin{equation}\label{lEaLq}
    \|S_{\beta}(t)u\|_{L^r} \leq Ct^{\frac{1}{\beta}\left(a+\frac{n}{r} \right)}\|u\|_{E},
    \end{equation}
for all $t>0$ and $u\in E$.
\end{prop}

\noindent
{\bf Proof.} Since $S_{\beta}(t_0):E \to L^{r_0}(\R^n)$ is bounded we have that $\|S_{\beta}(t_0)u\|_{L^{r_0}} \leq C\|u\|_{E}$ for every $u \in E$ and some constant $C>0$. As $L^{r_0}(\R^n)$ has scaling degree $-n/r_0$, by Proposition \ref{ers}, we have $\|S_{\beta}(t)u\|_{L^{r_0}}\leq Ct^{\frac{1}{\beta}\left(a+\frac{n}{r_0} \right)}\|u\|_{E}$ for all $t>0$. Furthermore, we know that $\|S_{\beta}(t)u\|_{L^r}\leq Ct^{-\frac{n}{\beta}\left( \frac{1}{r_0}-\frac{1}{r} \right)}\Vert u\Vert_{r_0},~1\leq r_0\leq r\leq +\infty$ and $t>0$. Thus
    \begin{eqnarray}
    \|S_{\beta}(t)u\|_{L^{r}} &\leq & C(t/2)^{-\frac{n}{\beta}\left( \frac{1}{r_0}-\frac{1}{r} \right)}\|S_{\beta}(t/2)u\|_{L^{r_0}} \nonumber\\
    &\leq & C'(t/2)^{-\frac{n}{\beta}\left( \frac{1}{r_0}-\frac{1}{r} \right)}(t/2)^{\frac{1}{\beta}\left(a+\frac{n}{r_0} \right)}\|u\|_{E} \nonumber\\
    &=& C''t^{\frac{1}{\beta}\left( a+\frac{n}{r} \right)}\|u\|_{E}, \nonumber
    \end{eqnarray}
for all $t>0$ and $u\in E$.

\begin{obs}\label{mathX} From here on, we denote by $\mathcal{X}$ the class of  Banach spaces  $E'$ which can be imbedded continuously in $\mathcal{S}'$ such that there exist $t_0$  and $r_0 \in [1,+\infty]$ and  $S_{\beta}(t_0): E \to L^{r_0}(\R^n)$ is bounded.
\end{obs}

\subsection{\bf The space $BE^{\alpha}$}\label{sessaoBEa} 
Let $\alpha \geq 0$ and let $E$ be a nontrivial Banach space which can be imbedded continuously in $\mathcal{S}'$ . We define $BE^{\alpha}$ as 
\begin{equation}\label{defnormaBEa}
 BE^\alpha=\{ u \in \mathcal{S}';    \|u\|_{BE^{\alpha}}= \sup_{t>0}t^{\alpha/\beta}\|S_{\beta}(t)u\|_E < \infty\}.
 \end{equation}
 It is clear that $BE^\alpha$ is a linear space and $\|\cdot\|_{BE^{\alpha}}$ defines a norm on $BE^{\alpha}$. Indeed, $\|u\|_{BE^{\alpha}}=0$ implies $ \|S_{\beta}(t)  u\|_{E}=0$, for every $t>0$. Since  $E$ is imbedded continuously in $\mathcal{S}'$ is continuous we conclude that $S_{\beta}(t)  u\to 0$ in $\mathcal{S}'$. Since  $S_{\beta}(t)  u\to u$ in $\mathcal{S}'$ as $t\to 0$, follows that $u=0.$ The other axioms are easily verified.

\begin{obs} 
\begin{enumerate}[(i)]
\item If $E$ is a Banach space with scaling degree $a$, then $\Vert \cdot \Vert_{BE^\alpha}$ has  scaling degree $a-\alpha$. Indeed, for $\lambda>0$ and Proposition \ref{sso}
$$
\begin{array}{ll}
\|u_{\lambda}\|_{BE^{\alpha}} &= \sup_{t>0}t^{\alpha/\beta}\|S_{\beta}(t) u_{\lambda}\|_{E}= \sup_{t>0}t^{\alpha/\beta}\|[S_{\beta}(\lambda^{\beta}t)u]_{\lambda}\|_{E} \\
&= \lambda^{a}\sup_{t>0}t^{\alpha/\beta}\|S_{\beta}(\lambda^{\beta}t)u\|_{E}\\
 &= \lambda^{(a-\alpha)}\sup_{t>0}(\lambda^{\beta}t)^{\alpha/\beta}\|S_{\beta}(\lambda^{\beta}t)u\|_{E}= \lambda^{(a-\alpha)}\|u\|_{BE^{\alpha}}.
 \end{array}
$$

\item If the norm $\Vert \cdot \Vert$ is translation invariant, then $\Vert \cdot \Vert_{BE^\alpha}$ is also one.  Indeed,  
$$
\begin{array}{ll}
 \|T_yu\|_{BE^{\alpha}} &= \sup_{t>0}t^{\alpha/\beta}\|S_{\beta}(t)T_yu\|_{E}=\sup_{t>0}t^{\alpha/\beta}\|T_{-y}S_{\beta}(t)u\|_{E} \\
&= \sup_{t>0}t^{\alpha/\beta}\|S_{\beta}(t)u\|_{E} = \|u\|_{BE^{\alpha}}. 
\end{array}
$$
\end{enumerate}
\end{obs}

Our objective now is to show that $BE^\alpha$ is a Banach space. To do this, we consider some properties of the homogeneous Besov space $\dot{B}^{\gamma}_{p,q}(\Rn)$ with $\gamma<0$. 

The following result was proved in \cite{miao}(Proposition 2.1).
\begin{prop}\label{lemamioabesov}
Let $1\leq r,s \leq +\infty$, $\gamma<0$ and $0<\beta<+\infty$.Then $u\in \dot{B}^{\gamma}_{r,s}(\Rn)$ if and only if
    \begin{equation}\nonumber
    \|u\|_{\dot{B}_{r,s}^{\gamma}} = \left\{
    \begin{array}{ll}
    \displaystyle \left[ \int_0^{+\infty}\left( t^{-\gamma/\beta}\|S_{\beta}(t)u\|_{L^r} \right)^s \frac{dt}{t}\right]^{\frac{1}{s}}, & \mbox{if }1\leq s<\infty, \\
    \displaystyle \sup_{t>0}t^{-\gamma/\beta}\|S_{\beta}(t)u\|_{L^r}, & \mbox{if }s=+\infty.
    \end{array}
    \right.
    \end{equation}
is finite.
\end{prop}
Note that if $E=L^r(\R^n)$ the set $BE^{\alpha}$, $\alpha>0$, is exactly the space $\dot{B}^{-\alpha}_{r,\infty}(\Rn)$.

In the next result we establish a continuous imbedding of  $BE^{\alpha}$ in a homogeneous Besov space.

\begin{prop}[\bf Imbedding in Besov spaces]\label{imersaoembesov} Let $E$ be a nontrivial Banach space continuously imbedded in $\mathcal{S}'$. Let $\alpha \geq 0$, $E \in \mathcal{X}$ with scaling degree $a$. There exist $r>1$ such that the imbedding $BE_a^{\alpha}\subset \dot{B}_{r,\infty}^{\gamma}(\R^n)$ is continuous with $\gamma = n/r+a-\alpha<0$.

\end{prop}

\noindent
{\bf Proof. } Since $E \in \mathcal{X}$,  let $t_0>0$ and $r_0\in [1,+\infty]$ so that $S_{\beta}(t_0):E\to L^{r_0}(\R^n)$ is bounded. From Proposition \ref{ers} we conclude that $a+n/r_0\leq 0$. Let $r>r_0$ so that $a+n/r<0$. From Proposition \ref{lemaqmq0}, there exists $C>0$ such that
    \begin{equation}\nonumber
    \|S_{\beta}(t)u\|_{L^r} = \|S_{\beta}(t/2)S_{\beta}(t/2)u\|_{L^r} \leq C(t/2)^{\frac{1}{\beta}(a+\frac{n}{r})}\|S_{\beta}(t/2)u\|_{E}.
    \end{equation}
Multiplying  this inequality by $\displaystyle t^{-\gamma/\beta}$ with $\gamma=n/r+a-\alpha$ we obtain
    \begin{equation}\nonumber
    t^{-\frac{\gamma}{\beta}}\|S_{\beta}(t)u\|_{L^{r}}\leq C2^{-\frac{\gamma}{\beta}}(t/2)^{\frac{\alpha}{\beta}}\|S_{\beta}(t/2)u\|_{E} = C'(t/2)^{\frac{\alpha}{\beta}}\|S_{\beta}(t/2)u\|_{E}.
    \end{equation}
Hence,  we get $\sup_{t>0}t^{-\gamma/\beta}\|S_{\beta}(t)  u\|_{L^{r}} \leq C'\|u\|_{BE^{\alpha}}$.

Now, the result follows from  Proposition \ref{lemamioabesov}.


Finally, we show that $BE^\alpha$ is a Banach space.
\begin{prop}[\bf Completeness  of $BE^{\alpha}$]\label{BEacompleto}
Let $E$ be a nontrivial Banach space continuously imbedeed  in $\mathcal{S}'$. Let $\alpha\geq 0$, $E \in \mathcal{X}$ with scaling degree $a$. Then the space $BE^{\alpha}$ is a Banach space. Furthermore, there exists $r>1$ so that the imbedding $ BE^\alpha \subset \dot{B}^{\gamma}_{r,\infty}(\Rn)$ is continuous, with $\gamma=n/r+a-\alpha$.
\end{prop}

\noindent
{\bf Proof. }Since  $E\subset \mathcal{S}'$ is a normed linear space it is suffices to show that $BE^{\alpha}$ is complete. Let $(u_n)_{n\in\N}$ a Cauchy sequence in $BE^\alpha$. For every $t>0$ the sequence $(S_{\beta}(t)  u_n)_{n\in \N}$ is a Cauchy sequence in $E$ because $$\|S_{\beta}(t)  u_n-S_{\beta}(t)  u_m\|_{E}\leq t^{-\frac{\alpha}{\beta}}\|u_n-u_m\|_{BE^{\alpha}}.$$ Hence, since $E$ is a Banach space we have $u(t):=\lim_{n\to\infty}S_{\beta}(t)  u_n$ in  $E$. Using the fact  that  the embedding $E \subset \mathcal{S}'$ is continuous,we conclude that $S_{\beta}(t)  u_n\to u(t)$ in $\mathcal{S}'$.

On the other hand, since $(u_n)_{n\in\N}$ is a Cauchy sequence in $BE^{\alpha}$ we obtain a constant $C>0$ such that $\|u_n(t)\|_{E}\leq Ct^{-\frac{\alpha}{\beta}}$ for every $t>0$. Thus,  $\|u(t)\|_{E}\leq Ct^{-\frac{\alpha}{\beta}}$ for every $ t>0$. 

We show now that there exists $v\in \mathcal{S}'$ so that $u(t)=S_{\beta}(t)  v$. From  Proposition \ref{imersaoembesov} the embedding $BE^{\alpha}\subset \dot{B}_{r,\infty}^{\gamma}$ is continuous for some $r>1$, $\gamma=n/r+a-\alpha$. Therefore, the Cauchy sequence $(u_n)_{n\in \N}$ converge in $\dot{B}_{r,\infty}^{\gamma}$ for some function $v$. Since  $\dot{B}_{r,\infty}^{\gamma}\subset \mathcal{S}'$ we have that $u_n\to u$ and $S_{\beta}(t)  u_n\to S_{\beta}(t)  v$ in $\mathcal{S}'$. By uniqueness we conclude that $u(t)=S_{\beta}(t)  v$.

\subsection{Adequate spaces} \label{Ade.eq} Let $E$ be a Banach space and $B$ a $p-$linear form defined on $E^p$. We say that $E$ is adequate to problem (\ref{In.uno}) if
\begin{enumerate}[(i)]
\item $\mathcal{S}\subset E \subset \mathcal{S}'$ both with continuous injections.
\item The norm $\|\cdot\|_E$ is invariant by translations, that is, $\|T_yu\|_E=\|u\|_E$ for every $u \in E$, $y \in \Rn$ and $T_yu=u(\cdot -y)$.
\item For every $u_i \in E, ~i=1,...,p$, we have $B(u_1,...,u_{p})\in \mathcal{S}'$. Moreover,
$$ \|S_{\beta}(t)  B(u_1,...,u_{p})\|_{E}\leq \omega(t)\prod_{i=1}^{p}\|u_i\|_{E}$$
where $\omega:(0,+\infty) \to (0,+\infty)$ and $\omega \in L^1(0,T)$ for  $0<T<+\infty$.
\end{enumerate}

\begin{obs}\label{sgee}If $~(E,\|\cdot\|_E)$ is a  Banach space satisfying conditions (i) and (ii), then
    \begin{equation}\label{Lim.s}
    \|S_{\beta}(t)u\|_{E} \leq \|K_{\beta}(1, \cdot)\|_{L^1}\|u\|_{E},
    \end{equation}
since  $K_{\beta}(1, \cdot) \in L^1(\R^n)$, see Lemma 2.1 of \cite{miao}). In particular, $E$ is adequate to problem (\ref{In.uno}) for $B(u,...,u)=u$.

\end{obs}

In the next result we establish estimates for $S_\beta(t) B(u_1,...,u_p)$ in the spaces $E$ and $BE^\alpha$.

\begin{prop}\label{esb1}
Let $B$ be a $p-$linear form with scaling degree $\sigma$ and let $E$ be a Banach space adequate to problem (\ref{In.uno}) with scaling degree $a$.
\begin{enumerate}[(i)]
\item  There exists a constant $C_1=C_1(a, \sigma,\beta, p)>0$ such that
    \begin{equation}\label{estsgea}
    \|S_{\beta}(t) B(u_1,...,u_{p})\|_{E}\leq C_1t^{[(p-1)a-\sigma]/\beta}\prod_{i=1}^{p}\|u_i\|_{E}.
    \end{equation}

\item Assume that  $0\leq \alpha \leq -(p-1)a + \sigma.$ Then there exists $C_2=C_2(a,\sigma,\beta, p)>0$ such that
    \begin{equation}\label{estsgbea}
    \|S_{\beta}(t) B(u_1,...,u_{p})\|_{BE^{\alpha}}\leq C_2t^{[\alpha +(p-1)a-\sigma]/\beta}\prod_{i=1}^{p}\|u_i\|_{E}.
    \end{equation}
\end{enumerate}
\end{prop}

\noindent{\bf Proof. } (i) Let $u_i\in E_a,~i=1,...,p$ and $t_0>0$. Since $B$ has scaling degree $\sigma$, we obtain from Proposition \ref{sso}, $S_{\beta}(t_0)B((u_1)_{\lambda},...,(u_{p})_{\lambda}) = \lambda^{\sigma}[S_{\beta}(\lambda^{\beta}t_0)B(u_1,...,u_{p})]_{\lambda}.$ Hence, using  Proposition \ref{sso} again,  the facts that $B$ is adequate and $E$ has scaling degree $a$, we have
$$
\begin{array}{ll}
    \lambda^{\sigma+a}\|S_{\beta}(\lambda^{\beta}t_0)B(u_1,...,u_{p})\|_{E} &=\|\lambda^{\sigma}[S_{\beta}(\lambda^{\beta}t_0)B(u_1,...,u_{p})]_{\lambda}\|_{E} \\
    &= \|S_{\beta}(t_0)B((u_{\lambda})_1,...,(u_{\lambda})_{p})\|_{E} \\ 
    &\leq \omega(t_0)\prod_{i=1}^{p}\|(u_i)_{\lambda}\|_{E} \\ 
    &= \omega(t_0)\lambda^{pa}\prod_{i=1}^{p}\|u_i\|_{E}.
\end{array}
$$
Thus, $\|S_{\beta}(\lambda^{\beta}t_0)B(u_1,...,u_{p})\|_{E}\leq \omega(t_0)\lambda^{[(p-1)a-\sigma]}\prod_{i=1}^{p}\|u_i\|_{E},$
for all $\lambda>0$. Setting $\lambda=(t/t_0)^{\frac{1}{\beta}}$ with $t>0$, we obtain (\ref{estsgea}) with $C_1=\omega(t_0)t_0^{\frac{-(p-1)a+\sigma}{\beta}}$.

(ii)  We first claim that there exists a constant $C'=C'(t)>0$ such that for all $\tau>0$,
    \begin{equation}\nonumber
    \tau^{\frac{\alpha}{\beta}}\|S_{\beta}(\tau)S_{\beta}(t)B(u_1,...,u_p)\|_{E}\leq C'(t)\prod_{i=1}^{p}\|u_i\|_{E}.
    \end{equation}
Indeed, since $\alpha\geq 0$, (\ref{Lim.s}) holds and $B$ is adequate, we have for  $0<\tau\leq 1$,
    \begin{eqnarray}\tau^{\frac{\alpha}{\beta}}\|S_{\beta}(\tau)S_{\beta}(t) B(u_1,...,u_{p})\|_{E_a} &\leq & \|S_{\beta}(\tau)S_{\beta}(t) B(u_1,...,u_{p})\|_{E_a}\nonumber \\ &\leq & \overline{C}\|S_{\beta}(t) B(u_1,...,u_{p})\|_{E_a} \nonumber \\ &\leq & \overline{C}\omega(t)\prod_{i=1}^{p}\|u_i\|_{E_a}. \nonumber
    \end{eqnarray}
On the other hand, from bound (\ref{Lim.s}), estimate (\ref{estsgea}) and $\alpha+(p-1)a-\sigma\leq 0$ we obtain for $\tau>1$, 
 $$
 \begin{array}{ll}
    \tau^{\frac{\alpha}{\beta}}\|S_{\beta}(\tau)S_{\beta}(t) B(u_1,...,u_{p})\|_{E} &= \tau^{\frac{\alpha}{\beta}}\|S_{\beta}(t) S_{\beta}(\tau) B(u_1,...,u_{p})\|_{E}\\
    &\leq  \til{C}\tau^{\frac{\alpha}{\beta}}\|S_{\beta}(\tau)B(u_1,...,u_{p})\|_{E}\\
    &\leq  \til{C}\tau^{\frac{\alpha+(p-1)a-\sigma}{\beta}}\prod_{i=1}^{p}\|u_i\|_{E}\\
    &\leq  \til{C}\prod_{i=1}^{p}\|u_i\|_{E}. 
    \end{array}
$$
Thus,  the claim holds for $C'(t)=\max\{\overline{C}w(t), \til{C}\}$. 

Fix now $t_0>0$. Since
    \begin{eqnarray}
    C'(t_0)\lambda^{pa}\prod_{i=1}^{p}\|u_i\|_{E} &=& C'(t_0)\prod_{i=1}^{p}\|(u_i)_{\lambda}\|_{E}\nonumber\\
    &\geq& \tau^{\frac{\alpha}{\beta}}\|S_{\beta}(\tau)S_{\beta}(t_0)B((u_{\lambda})_1,...,(u_{\lambda})_{p})\|_{E} \nonumber \\
    &=& \lambda^{\sigma+a}\tau^{\frac{\alpha}{\beta}}\|S_{\beta}(\lambda^{\beta}\tau)S_{\beta}(\lambda^{\beta}t_0)B(u_1,...,u_{p})\|_{E} \nonumber\\
    &=& \lambda^{\sigma+a-\alpha}(\lambda^{\beta}\tau)^{\frac{\alpha}{\beta}}\|S_{\beta}(\lambda^{\beta}\tau)S_{\beta}(\lambda^{\beta}t_0)B(u_1,...,u_{p})\|_{E}\nonumber
    \end{eqnarray}
we have $(\lambda^{\beta}\tau)^{\alpha/\beta}\|S_{\beta}(\lambda^{\beta}\tau)S_{\beta}(\lambda^{\beta}t_0)B(u_1,...,u_{p})\|_{E} \leq C'(t_0)\lambda^{\alpha+(p-1)a-\sigma}\prod_{i=1}^{p}\|u_i\|_{E}.$ Taking the supreme on $\tau$ and setting $\lambda=(t/t_0)^{1/\beta}$ we get (\ref{estsgbea}) with
$C_2=C'(t_0)t_0^{\frac{-(p-1)a+\sigma-\alpha}{\beta}}$.


\section{Local and global existence for problem (\ref{In.uno})}

The existence, local and global, of solutions for problem (\ref{In.uno}) is based in following abstract result.
\begin{lem}\label{labst}
Assume that $X$ is a Banach space and $A:X\times... \times X\to X$ is a $p-$linear form($p>1$) verifying
    \begin{equation}\label{dpa2}
    \|A(u_1,...,u_{p})\| \leq K\prod_{i=1}^{p}\|u_i\|,
    \end{equation}
for all $u_i\in X,~i=1,...,p$ and for some constant $K>0$. Let $M, R>0$ such that 
\begin{equation}\label{abs.dos}
R+pKM^{p}<M.
\end{equation}
 Then, for every $y\in X$ with $\|y\|\leq R$ the equation
    \begin{equation}\label{eqabst}
    u=y+A(u,...,u)
    \end{equation}
has a unique solution $u\in X$ and  $\|u\|\leq M$. Moreover, the solution $u$ depends continuously in the sense that, if  $v$ is a solution of (\ref{eqabst}), with $y_1$ in place of $y$, and $\Vert y_1\Vert\leq R$, $\Vert v\Vert \leq M$, then
    \begin{equation}\label{cabst}
    \|u-v\|\leq (1-pKM^{p-1})^{-1}\|y-y_1\|.
    \end{equation}
\end{lem}

\noindent{\bf Proof. } Set $B_{M} = \{u\in X;~\|u\| \leq M\}$. Consider the mapping $\mathcal{G}_y: B_{M} \to X$ defined by $ \mathcal{G}_y(u) = y + A(u,...,u).$
Since $A$ is $p$-linear and verify inequality (\ref{dpa2}) we deduce
$$ \|A(u,...,u)-A(v,...,v)\| \leq K \left( \sum_{k=0}^{p-1}\|u\|^{p-1-k}\|v\|^k \right) \|u-v\|.$$
Hence,
\begin{equation}\label{abs.uno}
\|\mathcal{G}_y(u)-\mathcal{G}_{y_1}(v)\| \leq \|y-y_1\| + K\left( \sum_{k=0}^{p-1}\|u\|^{p-1-k}\|v\|^k\right)\|u-v\|.
\end{equation}
Setting $y_1=v=0$ in (\ref{abs.uno}), we conclude by (\ref{abs.dos}) that $\Vert \mathcal{G}_y u\Vert \leq \Vert y\Vert+ K \Vert u\Vert^p\leq R+KM^{p}<M.$ From (\ref{abs.uno}), for $y_1=y$ we have that $\Vert \mathcal{G}_y(u)-\mathcal{G}_y(v)\Vert \leq p K M^{p-1}\Vert u -v\Vert$, where $KM^{p-1}<1$, by (\ref{abs.dos}). Therefore, $\mathcal{G}_y$ is a strict contraction on  $B_{M}$. Thus, $\mathcal{G}_y$ has a fixed point. The continuous dependence follows directly from (\ref{abs.uno}).

\begin{lem}\label{Lem.pre} Let $p>1$, $\beta>0$, $a, \sigma \in \R$ and $\alpha=(\beta-\sigma)/(p-1)+a.$ Assume that $0<\alpha+(p-1)a<\sigma$. Then,
\begin{enumerate}[(i)]
\item $\beta+(p-1)a=\sigma+(p-1)\alpha.$

\item $\beta+(p-1)a>\sigma,$ $\beta>p\alpha$ and $\beta+(p-1)a>\sigma-\alpha.$

\end{enumerate}
\end{lem}
\noindent{\bf Proof.} It follows directly since
$$
\begin{array}{rl}
\beta+(p-1)a-\sigma&=(p-1)\alpha>0\\
\beta-p\alpha&=-(p-1)a+\sigma-\alpha>0\\
\beta+(p-1)a -\sigma-\alpha&=p\alpha>0.
\end{array}
$$

\medskip \noindent 
{\bf Proof of the Theorem \ref{teorema1}.} Let  $ u_0 \in BE^{\alpha}$ and $$ X = L^{\infty}((0,\infty); BE^{\alpha}) \cap \left\{u:(0,\infty)\to E; \sup_{t>0}t^{\alpha/\beta} \|u(t)\|_{E} < \infty\right\}$$ with the norm 
$\|u\|_X= \sup_{t>0}\|u(t)\|_{BE^{\alpha}}+\sup_{t>0}t^{\alpha/\beta}\|u(t)\|_{E}.$
For $u_i \in E$, $i=1,2,...,p$ set
    \begin{equation}\label{Teo.uno}
    y=S_{\beta}(t) u_0\mbox{ and }
    A(u_1,...,u_p)(t)=\int_0^tS_{\beta}(t-\tau)B(u(\tau),...,u(\tau))d\tau.
    \end{equation}
Since $E \in \mathcal{X}$, by Proposition \ref{BEacompleto}, $BE^{\alpha}$ is a Banach space. Therefore,   $X$ is also a Banach space.  From Proposition \ref{esb1} and Lemma \ref{Lem.pre}
\begin{equation}\label{In.va}
\begin{array}{ll}
t^{\frac{\alpha}{\beta}}\Vert A(u_1,...,u_p)(t)\Vert_{E} &\leq  t^{\frac{\alpha}{\beta}}\int_0^t C_1(t-\tau)^{\frac{(p-1)a-\sigma}{\beta}}\prod_{i=1}^{p}\|u_i(\tau)\|_{E}d\tau \\
    &\leq  t^{\frac{\alpha}{\beta}}\left(\int_0^t C_1(t-\tau)^{\frac{(p-1)a-\sigma}{\beta}}\tau^{-\frac{p\alpha}{\beta}}d\tau\right)\prod_{i=1}^{p}\left(\sup_{t>0}t^{\frac{\alpha}{\beta}}\|u_i(t)\|_{E}\right) \\
    &= t^{1+\frac{(p-1)(a-\alpha)-\sigma}{\beta}}\left(\int_0^1 C_1(1-\tau)^{\frac{(p-1)a-\sigma}{\beta}}\tau^{-\frac{p\alpha}{\beta}}ds\right) \left(\prod_{i=1}^{p}\|u_i\|_{X}\right)\\
    &= K_1\prod_{i=1}^{p}\|u_i\|_{X},
    \end{array}
\end{equation}
where 
\begin{equation}\label{Val.k1}
K_1=\int_0^1 C_1(1-\tau)^{\frac{(p-1)a-\sigma}{\beta}}\tau^{-\frac{p\alpha}{\beta}}d\tau.
\end{equation}

Similarly, by Proposition  \ref{esb1} and the definition of $\alpha$ we obtain
\begin{equation}\label{In.tre}
\begin{array}{ll}
\|A(u_1,...,u_p)\|_{BE^{\alpha}} &\leq C_2 \int_0^t (t-\tau)^{\frac{(p-1)a-\sigma+\alpha}{\beta}}  \prod_{i=1}^{p}\|u_i(\tau)\|_{E}d\tau\\
&\leq C_2t^{1+\frac{(p-1)(a-\alpha)-\sigma}{\beta}}\prod_{i=1}^p\|u_i\||_X \int_0^1 (1-\tau)^{\frac{(p-1)a-\sigma+\alpha}{\beta}}\tau^{-\frac{p\alpha}{\beta}}d\tau\\
&=K_2 \prod_{i=1}^p\|u_i\|_X
\end{array}
\end{equation}
where 
\begin{equation}\label{Val.k2}
K_2=\int_0^1 C_2(1-\tau)^{-\frac{(p-1)a+\sigma-\alpha}{\beta}}\tau^{-\frac{p\alpha}{\beta}}d\tau.
\end{equation}
  Lemma \ref{Lem.pre} provides that  $K_1, K_2<\infty$. Hence, taking 
\begin{equation}\label{Val.k}
K=K_1+K_2
\end{equation} we conclude that $\|A(u_1,...,u_p)\|_{X}\leq K\prod_{i=1}^p\|u_i\|_{X}$.
From Lemma \ref{labst}, the  global existence and continuous dependence follows.

To show, the asymptotic behavior we argue as \cite{karch}. Arguing as (\ref{In.va}) and using (\ref{abs.uno}) it is possible to conclude 
$$
\begin{array}{ll}
t^{\alpha/\beta}\Vert u(t)-v(t)\Vert_E&\leq t^{\alpha/\beta}\Vert S_\beta(t)(u_0-v_0)\Vert_E+\\
& \hskip20pt C_1\int_0^t (t-\tau)^{\frac{(p-1)a-\sigma}{\beta}}\left ( \sum_{i=0}^p \Vert u \Vert^{p-1-k}_E \Vert v\vert^{k}_E\right)\Vert u(\tau)-v(\tau)\Vert_E d\tau\\
& \leq C_1pM^{p-1}\int_0^1(1-\tau)^{\frac{(p-1)a-\sigma}{\beta}}\tau^{-\frac{\alpha p}{\beta}}\left [(\tau t)^{\alpha/\beta}\Vert u(\tau t)-v(\tau t) \Vert_E\right]d\tau
\end{array}
$$
 Similarly, 
$$
\begin{array}{l}
t^{\alpha/\beta}\Vert u(t)-v(t)-S_\beta(t)(u_0-v_0)\Vert_E \\
\hskip20pt \leq C_1pM^{p-1}\int_0^1(1-\tau)^{\frac{(p-1)a-\sigma}{\beta}}\tau^{-\frac{\alpha p}{\beta}}\left [(\tau t)^{\alpha/\beta}\Vert u(\tau t)-v(\tau t) \Vert_E\right]d\tau.
\end{array}
$$
From these estimates and Lemma 6.1 of \cite{karch} the conclusion follows.

\noindent
{\bf Proof of Theorem \ref{teo.ex}.} Let  $ u_0 \in BE^{\alpha}$ and $$ X_T = L^{\infty}((0,T); BE^{\alpha}) \cap \left\{u:(0,T)\to E; \sup_{0<t<T}t^{\alpha/\beta} \|u(t)\|_{E} < \infty\right\}$$ with the norm 
$\|u\|= \sup_{0<t<T}\|u(t)\|_{BE^{\alpha}}+\sup_{0<t<T}t^{\alpha/\beta}\|u(t)\|_{E}.$
For $u_i \in E$, $i=1,2,...,p$ set $y$ and $A$ given by (\ref{Teo.uno}). Arguing as in the derivation of (\ref{In.va}) we obtain
$$t^{\alpha/\beta}\Vert A(u_1, ...,u_p)\Vert_E \leq K_1 T^{1+\frac{(p-1)(a-\alpha)-\sigma}{\beta}} \prod_{i=1}^p\|u_i\|_X,$$
where $K_1$ is given by (\ref{Val.k1}). Similarly, arguing as in the derivation of (\ref{In.tre}) we conclude
$$\Vert A(u_1, ...,u_p)\Vert_{BE^\alpha}\leq K_2 T^{1+\frac{(p-1)(a-\alpha)-\sigma}{\beta}} \prod_{i=1}^p\|u_i\|_X,$$ 
where $K_2$ is given by (\ref{Val.k2}). Hypotheses guarantee that constants $K_1,K_2$ are finite and that $1+[(p-1)(a-\alpha)-\sigma]/\beta>0$. Thus, we have
$\Vert A(u_1,...,u_p)\Vert_a \leq K \prod_{i=1}^p\|u_i\|_X,$
with $K=(K_1+K_2)T^{1+\frac{(p-1)(a-\alpha)-\sigma}{\beta}}.$ From Lemma \ref{labst} we have the desired result. 

\section{Global existence for system  (\ref{sistint1})}

 We extend the concept of adequate space for problem (\ref{In.uno}) given in subsection \ref{Ade.eq}. Let $E$ and $F$ be Banach spaces and let $B_1$ and  $B_2$ be $q-$linear form and $p-$linear form respectively.  We say that $E\times F$ is adequate to system (\ref{In.sis}), if 
\begin{enumerate}[(i)]
\item The inclusions  $\mathcal{S}\subset E, F \subset \mathcal{S}'$ are continuous.

\item The norms $\|\cdot\|_E$ and $\|\cdot\|_F$  are invariants for translations.

\item $B_1(v_1,...,v_{q}), B_2(u_1,...,u_{p}) \in \mathcal{S}',$ for every $ u_i \in E$, $i=1,..,p$ and  $ v_j\in F, j=1,...,q$, and
    \begin{eqnarray}
    \|S_{\beta}(t)  B_1(v_1,...,v_{q})\|_{E}\leq \omega_1(t)\prod_{i=1}^{q}\|v_i\|_{F}, \nonumber\\
    \|S_{\beta}(t)  B_2(u_1,...,u_{p})\|_{F}\leq \omega_2(t)\prod_{i=1}^{p}\|u_i\|_{E} \nonumber
    \end{eqnarray}
where $\omega_1,\omega_2:(0,+\infty)\to (0,+\infty)$  so that $~\omega_1,\omega_2 \in L^1(0,T)$, for every $0<T<+\infty$.
\end{enumerate}

If we consider $B_1(v_1,...,v_q)=u_1u_2\cdots u_p$ and  $B_2(u_1,...,u_p)=u_1u_2\cdots u_p$, then the spaces $L^r(\Rn) \times L^s(\Rn)$ and $L^{(r,r_1)}(\Rn)\times L^{(s,s_1)}(\Rn)$, with $s=r(q+1)/(p+1)>n(pq-1)/[\beta(p+1)]$, are adequate to system (\ref{In.sis}). This fact, follows from the following estimates
$$
\begin{array}{ll}
\|S_{\beta}(t)B_1(v_1,v_2,...,v_q)\|_r &\leq Ct^{-\frac{n}{\beta}\left( \frac{q}{s}-\frac{1}{r} \right)}\|v_1v_2\cdots v_q\|_{\frac{s}{q}} \\
&\leq Ct^{-\frac{n}{\beta}\left( \frac{q}{s}-\frac{1}{r}\right)}\prod_{i=1}^q\|v_i\|_s, \nonumber\\
\end{array}
$$
$$
\begin{array}{ll}
\|S_{\beta}(t)B_1(v_1,v_2,...,v_q)\|_{(r,r_1)} &\leq Ct^{-\frac{n}{\beta}\left( \frac{q}{s}-\frac{1}{r} \right)}\|v_1v_2\cdots v_q\|_{(\frac{s}{q},r_1)} \\
&\leq Ct^{-\frac{n}{\beta}\left( \frac{q}{s}-\frac{1}{r}\right)}\prod_{i=1}^q\|v_i\|_{(s,s_1)}.
\end{array}
$$

\begin{lem}\label{sm01} Let $B_1$ and $B_2$ be $q$ and  $p-$linear forms with scaling degree $\sigma_1$ and $\sigma_2$ respectively. Assume that $E$ and $F$ are Banach spaces with scaling degree $a$ and $b$ respectively and that $E\times F$ is adequate to system (\ref{In.sis}). Then, 
\begin{enumerate}[(i)]
\item There exist positive constants $C_1$ and $C_2$ such that
    \begin{eqnarray}\nonumber
    \|S_{\beta}(t) B_1(v_1,...,v_{q})\|_{E} &\leq & C_1 t^{(qb-a-\sigma_1)/\beta}\prod_{i=1}^{q}\|v_i\|_{F} \nonumber\\
    \|S_{\beta}(t) B_2(u_1,...,u_{p})\|_{E} &\leq & C_2 t^{(pa-b-\sigma_2)/\beta}\prod_{i=1}^{p}\|u_i\|_{E} \nonumber
    \end{eqnarray}

\item If $\alpha_1, \alpha_2>0$ and
$$
     \alpha_1 + qb \leq a+\sigma_1, \ \alpha_2 + pa \leq b+\sigma_2.
 $$
Then there exist positive constants $C_1$ and $C_2$ such that
    \begin{eqnarray}\nonumber
    \|S_{\beta}(t) B_1(v_1,...,v_{q})\|_{BE^{\alpha_1}} &\leq & C_1 t^{(\alpha_1-a+qb-\sigma_1)/\beta}\prod_{i=1}^{q}\|v_i\|_{F} \nonumber\\
    \|S_{\beta}(t) B_2(u_1,...,u_{p})\|_{BF^{\alpha_2}} &\leq & C_2 t^{(\alpha_2-b+pa-\sigma_2)]/\beta}\prod_{i=1}^{p}\|u_i\|_{E} \nonumber
    \end{eqnarray}
\end{enumerate}
\end{lem}

\noindent{\bf Proof. } The proof follows the same arguments used in the proof of  Proposition \ref{esb1}.

\medskip

We need also of the following technical result.
\begin{lem}\label{labst2} Let $E$ and $F$ be Banach spaces. Consider the system  
    \begin{equation}\label{sab1}
    \left\{
    \begin{array}{rcl}
    x &=& x_0 + B_1(y,...,y) \\
    y &=& y_0 + B_2(x,...,x)
    \end{array}
    \right.
    \end{equation}
where $B_1:F\times\cdots\times F \to E$ and $B_2:E\times\cdots\times E \to F$ are a $p-$linear  and a $q-$linear forms respectively. Assume that there exist $K_1,K_2>0$ such that
    \begin{eqnarray}
    \|A_1(y_1,...,y_q)\|_E\leq K_1\prod_{i=1}^{q}\|y_i\|_F \nonumber\\
    \|A_2(x_1,...,x_p)\|_F\leq K_2\prod_{i=1}^{p}\|x_i\|_E \nonumber
    \end{eqnarray}
Let $M,R>0$ verifying
    \begin{equation}\label{escRM}
    R + qM^{q}K_1 + pM^{p}K_2 < M.
    \end{equation}
Then, for every $(x_0,y_0)\in E\times F$ so that  $\|(x_0,y_0)\|_{E\times F} = \|x_0\|_E+\|y_0\|_F\leq R$, there exists a unique solution $(x,y)\in E\times F$ for the system (\ref{sab1}) such that $\|(x,y)\|_{E\times F}\leq M$. Moreover, if $\Vert (\bar{x}_0,\bar{y}_0)\Vert_{E \times F}\leq R$ and $(\bar{x},\bar{y})$ is the corresponding solution  of (\ref{sab1}) with $\Vert (\bar{x}, \bar{y})\Vert_{E \times F}\leq M$, then
    \begin{equation}\label{ab3}
    \|(x,y)-(\bar{x},\bar{y})\|_{E\times F} \leq [1-(qM^{q-1}K_1 + pM^{p-1}K_2)]^{-1}\|(x_0,y_0)-(\bar{x}_0,\bar{y}_0)\|_{E\times F}.
    \end{equation}
\end{lem}

{\bf Proof.} Set $B_M=\left\{ (x,y)\in E\times F ;~ \|(x,y)\|_{E\times F} = \|x\|_E+\|y\|_F \leq M\right\}$. Define $\mathcal{G}_{(x_0,y_0)}: B_{M} \to E\times F$ by $~\mathcal{G}_{(x_0,y_0)}(x,y)=\left( F_{x_0}(x,y), G_{y_0}(x,y) \right),$ with
    \begin{equation}\nonumber
    \left\{
    \begin{array}{rcl}
    F_{x_0}(x,y) &=& x_0 + A_1(y,...,y), \\
    G_{y_0}(x,y) &=& y_0 + A_2(x,...,x).
    \end{array}
    \right.
    \end{equation}
Since
    \begin{eqnarray}
    \|A_1(y,...,y)-A_1(\bar{y}...,\bar{y})\|_{E} &\leq & K_1\left(\sum_{k=0}^{q-1}\|y\|_{F}^{q-1-k}\|\bar{y}\|_{F}^k\right)\|y-\bar{y}\|_{F}, \nonumber\\
    \|A_2(x,...,x)-A_2(\bar{x}...,\bar{x})\|_{F} &\leq & K_2\left(\sum_{k=0}^{p-1}\|x\|_{E}^{p-1-k}\|\bar{x}\|_{E}^k\right)\|x-\bar{x}\|_{E}, \nonumber
    \end{eqnarray}
we obtain
    \begin{eqnarray}
    \|F_{x_0}(x,y)-F_{\bar{x}_0}(\bar{x},\bar{y})\|_{E} &\leq & \|x_0-\bar{x}_0\|_E + qM^{q-1}K_1\|y-\bar{y}\|_{F}, \nonumber\\
    \|G_{y_0}(x,y)-G_{\bar{y}_0}(\bar{x},\bar{y})\|_{F} &\leq & \|y_0-\bar{y}_0\|_F + pM^{p-1}K_2\|x-\bar{x}\|_{E}. \nonumber
    \end{eqnarray}
Hence,
\begin{equation}\label{cont.uno}
    \begin{array}{l}
    \|\mathcal{G}_{(x_0,y_0)}(x,y)-\mathcal{G}_{(\bar{x}_0,\bar{y}_0)}(\bar{x},\bar{y})\|_{E\times F}\\ 
\leq \|(x_0,y_0)-(\bar{x}_0,\bar{y}_0)\|_{E\times F} + ( qM^{q-1}K_1 + pM^{p-1}K_2 )\|(x,y)-(\bar{x},\bar{y})\|_{E\times F}.
    \end{array}
\end{equation}
Thus, for  $(x_0,y_0)=(\bar{x}_0,\bar{y}_0)$ we have
    \begin{equation}\nonumber
    \|\mathcal{G}_{(x_0,y_0)}(x,y)-\mathcal{G}_{(x_0,y_0)}(\bar{x},\bar{y})\|_{E\times F} \leq \left[ qM^{q-1}K_1 + pM^{p-1}K_2 \right]\|(x,y)-(\bar{x},\bar{y})\|_{E\times F}.
    \end{equation}
Therefore, from inequality (\ref{escRM}), $\mathcal{G}_{(x_0,y_0)}$ is a strict contraction. 

On the other hand, for $(\bar{x}_0,\bar{y}_0)=(\bar{x},\bar{y})=(0,0)$
    \begin{equation}\nonumber
\begin{array}{ll}
    \|\mathcal{G}_{(x_0,y_0)}(x,y)\|_{E\times F} &\leq \|(x_0,y_0)\|_{E\times F} + \left[ qM^{q-1}K_1 + pM^{p-1}K_2 \right]\|(x,y)\|_{E\times F} \\
&\leq R+qM^{q}K_1 + pM^{p}K_2 \leq M.
\end{array}
    \end{equation}
So, $\mathcal{G}_{(x_0,y_0)}(B_M)\subseteq B_M$ and the existence  follows by fixed point theorem.

 Continuous dependence follows from (\ref{cont.uno}).

\begin{lem}\label{Lem.ks} Let $p, q \geq 1$ such that $pq>1$, $\beta>0, a, b, \sigma_1, \sigma _2 \in \R$. Set 
$$ \alpha_1=\frac{\beta(q+1)}{pq-1}+a -\frac{\sigma_1+q\sigma_2}{pq-1}, \ \alpha_2=\frac{\beta(p+1)}{pq-1}+b-\frac{\sigma_2+p\sigma_1}{pq-1}.$$ 
Suppose that 
\begin{enumerate}[(i)]
\item $\alpha_1, \alpha_2>0$.

\item $\alpha_1+qb<a+\sigma_1$, $\alpha_2+pa<b+\sigma_2$.

\item $\alpha_1<q\alpha_2$, $\alpha_2<p\alpha_1$.
\end{enumerate}

Then,
\begin{enumerate}[(i)]
\item $\beta-a+ \alpha_1-\alpha_2q=-q b+\sigma_1,$ $\beta-b-\alpha_1 p+\alpha_2=-p a+\sigma_2.$

\item $\beta-a>-qb + \sigma_1,$ $\beta-b>-pa+\sigma_2.$

\item $\beta>q \alpha_2$, $\beta>p\alpha_1$.

\item $\beta+\alpha_1-a>-q b+\sigma_1$, $\beta+\alpha_2-b>-pa+\sigma_2.$
\end{enumerate}
\end{lem}
\noindent{\bf Proof. } Follows directly since
$$
\begin{array}{rl}
\beta-a+qb-\sigma_1&=q \alpha_2-\alpha_1\\
\beta-q\alpha_2&=-\alpha_1+a-qb+\sigma_1\\
\beta +\alpha_1-a+qb-\sigma_1&=q\alpha_2.
\end{array}
$$
Similar argument can be used tho show the other inequalities.  

\bigskip
{\bf Proof fo Theorem \ref{teo2}} Let $\alpha_1,\alpha_2$ defined by (\ref{defalpha12}), and let $(u_0,v_0) \in BE^{\alpha_1}\times BF^{\alpha_2}$. Set $x_0=S_{\beta}(t)u_0$, $y_0=S_{\beta}(t)v_0$ and 
    \begin{eqnarray}
    A_1(v_1,...,v_q)(t) &=& \int_0^t S_{\beta}(t-\tau)B_1(v_1(\tau),...,v_q(\tau))d\tau, \nonumber\\
    A_2(u_1,...,u_p)(t) &=& \int_0^t S_{\beta}(t-\tau)B_2(u_1(\tau),...,u_p(\tau))d\tau. \nonumber
    \end{eqnarray}
 Consider $X, Y$ given by
    \begin{eqnarray}
    X &=& \left\{u: (0,+\infty) \to E;  \|u\|_X=\sup_{t>0}t^{\alpha_1/\beta}\|u(t)\|_{E} + \sup_{t>0}\|u(t)\|_{BE^{\alpha_1}} < +\infty \right\}, \nonumber\\
    Y &=& \left\{v: (0,+\infty) \to F; \|u\|_2=\sup_{t>0}t^{\alpha_2/\beta}\|v(t)\|_{F} + \sup_{t>0}\|v(t)\|_{BF^{\alpha_2}}  < +\infty\right\}. \nonumber
    \end{eqnarray}
Since $E, F \in \mathcal{X}$, from Proposition \ref{BEacompleto}, we conclude that $BE^{\alpha_1}$ and $BF^{\alpha_2}$ are Banach spaces. Therefore,  $(X,\|\cdot\|_1)$ and $(Y,\|\cdot\|_2)$ are Banach spaces. Thus,  $X\times Y$ is also a Banach space with the norm
    \begin{equation}\nonumber
    \|(u,v)\|_{X\times Y}:=\|u\|_X+\|v\|_Y.
    \end{equation}
Since $\alpha_1 + qb<a+\sigma_1$ and $\alpha_2 + pa \leq b+ \sigma_2$, Lemma \ref{sm01} can be used. Note that
$\beta+\alpha_1+qb-a -\sigma_1-\alpha_2 q=0,$
$\beta+\alpha_2+pa-b-\sigma_2-p\alpha_1=0.$
By Lemma \ref{sm01}(i) we have
    \begin{eqnarray}
    t^{\alpha_1/\beta}\|A_1(v_1,...,v_q)(t)\|_{E} &\leq &  K_1 \left(\prod_{i=1}^q\|v_i\|_2\right), \label{ak1}\\
    t^{\alpha_2/\beta}\|A_2(u_1,...,u_p)(t)\|_{F} &\leq &  K_2 \left(\prod_{i=1}^p\|u_i\|_1\right), \label{ak2}
    \end{eqnarray}
where
    \begin{eqnarray}
    K_1 &=& C_1\int_0^1(1-s)^{\frac{qb-a-\sigma_1}{\beta}}s^{-\frac{q\alpha_2}{\beta}}ds < +\infty,\nonumber\\
    K_2 &=& C_2\int_0^1(1-s)^{\frac{pa-b-\sigma_2}{\beta}}s^{-\frac{p\alpha_1}{\beta}}ds < +\infty.\nonumber
    \end{eqnarray}
By Lemma \ref{sm01}(ii) we have
    \begin{eqnarray}\label{estimativa1}
    \|A_1(v_1,...,v_q)(t)\|_{BE^{\alpha_1}} 
    &\leq & \til{K}_1\left(\prod_{i=1}^q\|v_i\|_2\right), \label{ak1til}\\
    \|A_2(u_1,...,u_p)(t)\|_{BF^{\alpha_2}} 
    &\leq & \til{K}_2\left(\prod_{i=1}^p\|u_i\|_{1}\right), \label{ak2til}
    \end{eqnarray}
where
    \begin{eqnarray}
    \til{K}_1=C_1\int_0^1 (1-s)^{\frac{\alpha_1+qb-a-\sigma_1}{\beta}}s^{-\frac{q\alpha_2}{\beta}}ds, \nonumber\\
    \til{K}_2=C_2\int_0^1 (1-s)^{\frac{\alpha_2+pa-b-\sigma_2}{\beta}}s^{-\frac{p\alpha_1}{\beta}}ds, \nonumber
    \end{eqnarray}
Therefore,
    \begin{eqnarray}
    \|A_1(v_1,...,v_q)\|_1\leq K'_1\left( \prod_{i=1}^{p}\|v_i\|_2 \right), \nonumber\\
    \|A_2(u_1,...,u_p)\|_2\leq K'_2\left( \prod_{i=1}^{p}\|u_i\|_1 \right), \nonumber
    \end{eqnarray}
where  $K'_i=K_i+\til{K}_i,~i=1,2.$ A finiteness  of $K_1,K_2,\til{K}_1,\til{K}_2$ are consequences of Lemma \ref{Lem.ks}.  Moreover, the right side of estimates (\ref{ak1})-(\ref{ak2til})  do not depend of $t$. Thus, the result follows from Lemma \ref{labst2}.


\begin{thebibliography}{99}


\bibitem{AEZ} J. Aguirre, M. Escobedo and E. Zuazua, Self-similar solutions of a convection diffusion equations and related elliptic problems, Comm. Partial Differential Equations 15, 139-157(1990). 

\bibitem{Ben} M. Ben-Artzi, P. Souplet and F. B. Weissler. The local theory for viscous Hamilton Jacobi equations in Lebesgue space, J. Math. Pures Appl. (9) 81, 343-378, (2002).

\bibitem{bergh} { J. Bergh, J. L\"ofstrom}, { Interpolation Spaces}, an introduction, Springer-Verlag, New York, 1976.

\bibitem{BC} H. Brezis and T. Cazenave, A nonlinear heat equation with singular initial data, J. Anal. Math. 68 (1996), 277-304.



\bibitem{Can} M. Cannone and G. Karch, About the regularized Navier-Stokes equations, J. Math. Fluid Mech. 7 (2005) 1-28.

\bibitem{cazenaveweissler} { T. Cazenave and  F. B. Weissler}, { Asymptotically self-similar global solutions of the nonlinear Schr\"odinger and heat equations}, Math. Z. 228, 83-120 (1998).


\bibitem{chandrasekhar} A.  S. Chandrasekhar, { Stochastic problems in physics and astronomy,} Rev. Mod. Phys. vol.{\bf 15} (1943), 1-89.

\bibitem{escobedoherrero1} { M. Escobedo, M.A. Herrero}, { Boundedness and blow up for a semilinear reaction diffusion system}, J. Differential Equations, 89 (1991), 176-202.




\bibitem{EZ} M. Escobedo and E. Zuazua, Large time behavior for convection-diffusion equations in $\R^N$, J. Differential Equations 100, 119-161(1991).

\bibitem{fujita} { H. Fujita}, { On the blowing up of solutions of the Cauchy problem for $u_t=\Delta u + u^{1+\alpha}$}, J. Fac. Sci. Univ. Tokyo Sect. A Math. 16 (1966) 105-113.

\bibitem{grafakos} { L. Grafakos}, { Modern Fourier Analysis}, Ed. 2a, Springer-Verlag, New York, 2009.

\bibitem{karch} { G. Karch}, { Scaling in nonlinear parabolic equations}, J. Math. Anal. Appl. 234, 534-558 (1999), 

\bibitem{karch2} { G. Karch}, { Scalling in nonlinear parabolic equations: locality versus globality}

\bibitem{kilbas} { A. A. Kilbas, H.M. Srivastava, J.J. Trujillo}, { Theory and Applications of Fractional Diferential Equations}, Vol. 204 (North-Holland Mathematics Studies), 2006.

\bibitem{wu} { C. Miao, B. Yuan}, { Solutions to some nonlinear parabolic equations in pseudomeasure spaces}, Math. Nachr. 280 (2007) 171-186.

\bibitem{miao} { C. Miao, B. Yuan, B. Zhang}, { Well-posedness of the Cauchy problem for the fractional power dissipative equations}, Nonlinear Anal. 68, 461-484 (2008)















\bibitem{roberts} { C. A. Roberts, W. E. Olmstead}, { Blow-up in a subdiffusive medium of infinite extent.} Fract. Calc. Appl. Anal. {\bf 12} (2009), no. 2, 179-194.

\bibitem{snoussi1} { S. Snoussi, S. Tayachi, Fred B. Weissler}, { Asymptotically self-similar global solutions of a semilinear parabolic equation with a nonlinear gradient term}, Proc. Royal Soc. Edinburgh Sect. A., 129 (1999), 419-440.

\bibitem{snoussi2} { S. Snoussi, S. Tayachi}, { Global existence, asymptotic behavior and self-similar solutions for a class of semilinear parabolic systems}, Nonlinear Analysis, 48 (2002), 13-35.


\bibitem{snoussi4}  S. Snoussi, S. Tayachi, Fred B. Weissler,  Asymptotically self-similar global solutions of a general semilinear heat equation, Math. Ann. 321,(2001) 131-155. 


\bibitem{wei1} F. B. Weissler, semilinear evolution equations in Banach spaces, J. Funct. Anal. 32 (1979), 277-296.


\bibitem{wei2} F. B. Weissler, Local existence and nonexistence for semilinear parabolic equations in $L^p$, Indiana Univ. Math. J. 29 (1980), 79-102. 

\bibitem{wuyuan} { G. Wu, Jia Yuan}, { Well-posedness of the Cauchy problem for the fractional power dissipative equation in critical Besov spaces}, J. Math. Anal. Appl. 340, 1326-1335 (2008).

\end{thebibliography}
\end{document}